\DocumentMetadata{
  lang=en-US,
  pdfstandard=UA-2,
  tagging=on,
  tagging-setup={
    math/setup={mathml-AF,mathml-SE},
    table/header-rows=1,
  extra-modules={verbatim-mo}}
}
\documentclass[12pt,letterpaper]{amsart}
\usepackage[american]{babel}
\usepackage[margin=1.0in, portrait]{geometry}

\usepackage{meyermacros}
\usepackage{tikz-cd}
\tikzcdset{ampersand replacement={\&}}

\usepackage[stixtwo]{fontsetup} 

\usepackage{tagpdf}
\tagpdfsetup{activate/all, math/alt/use}

\usepackage{graphicx}

\usepackage{array, fancyhdr}
\usepackage[dvipsnames]{xcolor}

\usepackage{ragged2e}
\title[Obstructing Aspherical $3$-Manifold Groups via Torus
Surgery]{An Obstruction to Realizing Aspherical $3$-Manifold Groups
via Torus Surgery}
\author{Nicholas Meyer}
\address{Department of Mathematics and Statistics\\ \mbox{ Winona
State University }\\ Winona, MN}
\email{nick.meyer@winona.edu}
\urladdr{https://nickmeyer.phd}
\date{November 13, 2026}

\usepackage{hyperref}

\usepackage[english=american]{csquotes}
\usepackage[backend=biber,maxcitenames=5,maxbibnames=9,doi=false,isbn=false,url=false,eprint=false,style=numeric]{biblatex}
\renewbibmacro{in:}{}
\hypersetup{
    allbordercolors=blue,
    pdfdisplaydoctitle=true,
    colorlinks=true,
    citecolor=blue, filecolor=blue, 
    linkcolor=blue, urlcolor=blue,
    pdfauthor={Nicholas Meyer},
    pdftitle={An Obstruction to Realizing Aspherical 3-Manifold Groups via Torus Surgery},
    pdfkeywords={{fibration}, {torus surgery}, {aspherical manifolds}},
    pdfsubject={{Mathematics/Topology}},
}

\begin{document}

\begin{abstract}
  Let $Y$ be a closed, connected, oriented, aspherical $3$-manifold
  with fundamental group $\pi = \pi_1(Y)$ and positive first Betti
  number $\beta_1(Y) > 0$, and let $X_0 = Y \times S^1$. We prove
  that no torus surgery on $X_0$ along a torus-fibered $T^2$-knot $T
  \subseteq X_0$ can yield a $4$-manifold $X = X_{0; T, f}$ with
  fundamental group $\pi$. The obstruction combines the vanishing of
  the Euler characteristic of $X$ with a classifying-map argument
  showing that the unimodular intersection form on $H^2(X; \Z)$ must
  admit an isotropic subgroup of rank $\beta_1(\pi)$, forcing
  $\chi(X) \geq 2$. Along the way, we also present an independent
  geometric proof via covering space theory, the Sphere Theorem, and
  Scott's Core Theorem that the fundamental group of a closed
  aspherical $3$-manifold is torsion-free.
\end{abstract}

\maketitle

\section{Introduction}\label{sec:intro}

Surgery along embedded surfaces is a fundamental tool in $4$-manifold
topology for constructing new smooth $4$-manifolds and investigating
the interplay between smooth structures, fibrations, and algebraic
invariants. In \emph{torus surgery} (or logarithmic transformation),
one modifies a smooth, oriented $4$-manifold $X_0$ by removing the
interior of a tubular neighborhood $\nu T \cong T^2 \times D^2$ of a
smoothly embedded $2$-torus $T$ of self-intersection zero and
regluing $T^2 \times D^2$ via an orientation-preserving
diffeomorphism \mbox{$f: \partial(T^2 \times D^2)$} $\to$ \mbox{$\partial(X_0
\setminus \mathring{\nu}T)$} of the $3$-torus boundary, producing the
surgered $4$-manifold
\[
  X_{0; T, f} = (X_0 \setminus \mathring{\nu}T) \cup_f (T^2 \times D^2).
\]
This four-dimensional analog of Dehn surgery has long played a
central role in constructing exotic smooth structures, beginning
with the work of Donaldson \autocite{Donaldson87} and Dolgachev
\autocite{Dolgachev2010} on elliptic surfaces, and is closely
related to \hbox{Fintushel--Stern} knot surgery \autocite{FS98} and its many
variations \autocites{Akbulut2002,AEMS2008,Yun2008,Ni17,Kim19}.

A recurring theme in the study of torus surgery is understanding how
surgery interacts with fibrations and the fundamental group of the
ambient $4$-manifold. In the setting of $S^1 \times S^3$, Zentner
\autocite{Zentner} proved that torus surgery along $S^1$ times the
unknot $\calU$ preserves the diffeomorphism type whenever the
resulting manifold is a homology $S^1 \times S^3$. In the author's
doctoral dissertation \autocite{meyerthesis} and subsequent work
\autocite{torus_fibrations}, this perspective was broadened by
establishing fibration extension theorems for $4$-manifolds with
$T^3$ boundary. In particular, theorems related to products $S^1
\times E_Y(K)$ of a circle with a fibered knot exterior in a closed, oriented
$3$-manifold $Y$, as well as $\Sigma_{g,1}$-bundles over $T^2$ whose
monodromy fixes $\partial\Sigma_{g,1}$ pointwise were proved.
In particular, \autocite{torus_fibrations} introduced the notion of a
\emph{torus-fibered $T^2$-knot}: a smoothly embedded $2$-torus $T
\into X_0$ of self-intersection zero whose exterior $X_0 \setminus
\mathring{\nu}T$ fibers over $T^2$ with monodromy fixing the boundary
of the fiber pointwise. Archetypal examples arise in product
$4$-manifolds $X_0 = Y \times S^1$ by taking $T = K \times S^1$ for a
fibered knot $K \subseteq Y$.

As shown in \autocites{meyerthesis,torus_fibrations}, performing
torus surgery along a torus-fibered $T^2$-knot $T \subseteq X_0$ always
yields a closed $4$-manifold $X = X_{0; T, f}$ that fibers over
$S^1$. Moreover, the bundle structure on the exterior imposes strict
topological constraints: any closed $4$-manifold admitting a
torus-fibered $T^2$-knot must have vanishing Euler characteristic and
signature.

In light of these results --- and considering that in $S^1 \times
S^3$ one can perform non-trivial torus surgeries along $S^1 \times
\calU$ to obtain $S^1$ times lens spaces \autocite{torus_fibrations},
or even recover $\pi_1 \cong \Z$ in the homology $S^1 \times S^3$
case \autocite{Zentner} --- it is natural to ask which fundamental
groups can be realized by torus surgery along a torus-fibered
$T^2$-knot in $X_0 = Y \times S^1$. In particular, when $Y$ is a
closed, oriented, aspherical $3$-manifold (so that $Y$ is a $K(\pi,
  1)$ space for $\pi = \pi_1(Y)$, as is the case for hyperbolic
$3$-manifolds, torus bundles, or nilmanifolds), can torus surgery on
$X_0 = Y \times S^1$ along a torus-fibered $T^2$-knot reduce the
fundamental group from $\pi \times \Z$ to $\pi$?

In this note, we show that whenever $\beta_1(Y) > 0$, there is a
cohomological obstruction preventing any such surgery from realizing
the fundamental group $\pi = \pi_1(Y)$.

\begin{theorem}[store=mainthm]\label{thm:main}
  Let $\pi = \pi_1(Y)$ be the fundamental group of some closed,
  connected, oriented, aspherical $3$-manifold $Y$ with positive
  first Betti number, and let $X_0 = Y \times S^1$.

  There is not a pair of a torus-fibered $T^2$-knot $T$ in $X_0$ and
  a diffeomorphism $f: T^3 \to T^3$ such that $X= X_{0; T, f}$ has
  fundamental group $\pi$.
\end{theorem}

The proof of Theorem~\ref{thm:main} combines the bundle-theoretic
Euler characteristic identity $\chi(X) = 0$ from
\autocite{torus_fibrations} with a classifying-map argument
exploiting the asphericity of $Y$. Specifically, if $\pi_1(X) \cong
\pi$, the classifying map $\kappa: X \to Y$ induces an injection
\mbox{$\kappa^\ast: H^2(Y; \Z) \into H^2(X; \Z)$} whose image $L$ is an
isotropic subgroup of rank $\beta_1(\pi) \geq 1$ with respect to the
unimodular intersection form $Q_X$. Nondegeneracy of $Q_X$ then
forces \mbox{$\beta_2(X) \geq 2\beta_1(\pi)$,} which implies $\chi(X) \geq
2$, contradicting $\chi(X) = 0$.

The paper is organized as follows. In Section~\ref{sec:alg}, we
record an algebraic preliminary that the fundamental group of a
closed aspherical $3$-manifold is torsion-free, giving an independent
geometric proof via covering space theory, the Sphere Theorem
\autocite{Papa57}, and Scott's Core Theorem \autocite{Scott73}. In
Section~\ref{sec:main_thm}, we prove Theorem~\ref{thm:main}. In
Section~\ref{sec:examples}, we illustrate the obstruction across
several families of $3$- and $4$-manifolds and show that the
asphericity hypothesis on $Y$ is sharp, and in
Section~\ref{sec:conclusion}, we discuss concluding remarks and
further consequences.

\section{Algebraic Preliminaries}\label{sec:alg}

Before turning to four-dimensional torus surgery in
Section~\ref{sec:main_thm}, we record a foundational property of the
fundamental groups of closed aspherical $3$-manifolds. Classically,
the fact that the fundamental group $G$ of any finite-dimensional
aspherical CW complex $K(G, 1)$ is torsion-free is deduced
algebraically via group homology and cellular chain complexes: if $G$
contained a non-trivial finite cyclic subgroup $H \cong \Z_p$, the
corresponding covering space would be a finite-dimensional $K(H, 1)$,
contradicting the fact that $\Z_p$ has non-vanishing group homology
in infinitely many degrees (see, for instance,
\autocite[Example~2.43]{hatcher}).

When the aspherical space in question is a closed $3$-manifold,
however, one can bypass homological algebra altogether and establish
torsion-freeness directly from three-dimensional geometric topology.
As this argument was found independently and highlights the interplay
between covering spaces, irreducibility, and compact cores, we
include this alternate geometric proof here.
\clearpage
\begin{lemma}\label{lem:asph}
  If $G$ is an aspherical $3$-manifold group, then it is torsion-free.
\end{lemma}
\begin{proof}
  Let $M$ be a closed, orientable, aspherical $3$-manifold with
  fundamental group \mbox{$G = \pi_1(M)$,} so $\pi_k(M) = 0$ for all $k \ge
  2$ (if $M$ is non-orientable, we may pass to its orientation double
  cover, whose fundamental group is an index-$2$ subgroup of $G$). By
  Papakyriakopoulos's Sphere Theorem \autocite{Papa57}, an orientable
  3-manifold $Y$ satisfies $\pi_2(Y) = 0$ if and only if it is
  irreducible. The universal cover $\widetilde{M}$ is contractible
  (by Whitehead's Theorem \autocites{Whitehead1, Whitehead2}) and has
  empty boundary (by direct construction). Since a non-empty compact
  manifold (without boundary) cannot be contractible, $\widetilde{M}$
  must be non-compact, which forces $G$ to be infinite.

  Suppose $G$ contains a non-trivial element of finite order $p$,
  generating a finite cyclic subgroup \mbox{$H = \gen{g} \cong
  \Z_p$.} Consider the covering space $M_H = \widetilde{M} / H$.

  Because the index $\br{G:H}$ is infinite, $M_H$ is a connected,
  open 3-manifold without boundary, with finite fundamental group
  $\pi_1(M_H) \cong \Z_p$. Furthermore, covering space theory gives
  \mbox{$\pi_2(M_H) \cong \pi_2(M) = 0$,} hence $M_H$ is irreducible by the
  Sphere Theorem.

  By Scott's Core Theorem \autocite{Scott73}, every 3-manifold with a
  finitely generated fundamental group contains a compact core: a
  compact 3-manifold subdomain $C \subseteq M_H$ such that the
  inclusion map induces an isomorphism \mbox{$\pi_1(C) \to \pi_1(M_H)
  \cong \Z_p$.} Because $M_H$ is non-compact, the boundary $\partial
  C$ must be non-empty. If $\partial C$ contains any surface
  component $\Sigma$ of genus $g \ge 1$, its fundamental group
  $\pi_1(\Sigma)$ is infinite and torsion-free. From the commutative
  diagram of inclusions and covering maps in Figure~\ref{fig:cd1},
  \begin{figure}[!h]
    \centering\[
      \begin{tikzcd}[column sep=large]
        \Sigma \arrow[r, hook] \arrow[d, "\pi_1"']
        \& C \arrow[r, hook] \arrow[d, "\pi_1"']
        \& M_H = \widetilde{M}/H \arrow[r, "\mathrm{cov}"] \arrow[d, "\pi_1"']
        \& M \arrow[d, "\pi_1"] \\
        \pi_1(\Sigma) \arrow[r]
        \& \pi_1(C) \arrow[r, "\cong"]
        \& H \cong \Z_p \arrow[r, hook]
        \& G
    \end{tikzcd}\]
    \caption{The commutative diagram of inclusions and covering maps
    as outlined in this proof.}\label{fig:cd1}
  \end{figure}
  the induced map \mbox{$\pi_1(\Sigma) \to \pi_1(C)$} cannot be injective as
  $\pi_1(C)$ is finite. By the Loop Theorem and Dehn's Lemma
  \autocite{Papa57}, there exists a simple closed curve on $\Sigma$
  that is non-trivial in $\pi_1(\Sigma)$ but bounds an embedded
  2-disk in $C$. Compressing $\Sigma$ along a maximal disjoint
  collection of such disks simplifies the boundary without creating
  new non-trivial fundamental group elements. Carrying this out
  across all higher-genus boundary components replaces $C$ with a
  compact core $C' \subseteq M_H$ whose boundary $\partial C'$
  consists entirely of embedded 2-spheres.

  Because $M_H$ is irreducible, the Sphere Theorem guarantees that
  every 2-sphere boundary component of $\partial C'$ bounds an
  embedded 3-ball in $M_H$. These 3-balls cannot lie inside $C'$, as
  filling an interior void would simply eliminate a disconnected
  cavity rather than a boundary component bounding the non-compact
  exterior. Therefore, each boundary sphere of $\partial C'$ bounds
  an exterior 3-ball in $M_H \setminus \mathring{C'}$. Capping off
  every boundary component of $C'$ with its corresponding exterior
  3-ball constructs a compact 3-manifold with empty boundary inside
  $M_H$. Since $M_H$ is connected, this compact submanifold must be
  all of $M_H$, proving that $M_H$ is compact. This directly
  contradicts $\br{G:H} = \infty$. Thus, $G$ has no element of finite
  order, so $G$ is torsion-free.
\end{proof}

Beyond guaranteeing that $\pi = \pi_1(Y)$ is infinite and
torsion-free, the asphericity of $Y$ means that $Y$ is an
Eilenberg--MacLane space $K(\pi, 1)$ and therefore serves as a
geometric model for the classifying space $B\pi$. In particular, any
manifold $X$ with fundamental group $\pi$ admits a canonical homotopy
class of classifying maps $X \to Y$, allowing us to pull
cohomological information directly from $Y$ into $X$.

\section{Proof of Main Theorem}\label{sec:main_thm}

We are now ready to establish our main obstruction to realizing
aspherical $3$-manifold groups with positive first Betti number via
torus surgery along torus-fibered $T^2$-knots. The key tension lies
between the bundle structure of the knot exterior --- which forces
the surgered $4$-manifold $X = X_{0; T, f}$ to have vanishing Euler
characteristic --- and the unimodularity of the intersection form
$Q_X$ in the presence of the isotropic cohomology classes pulled back
from $Y$. For convenience, we restate Theorem~\ref{thm:main} before
giving its proof.

\getkeytheorem{mainthm}

\begin{proof}
  We do this by contradiction. Assume, for the sake of contradiction,
  that $T$ is a torus-fibered $T^2$-knot in $X_0$ and $f: T^3 \to
  T^3$ is an orientation-preserving diffeomorphism such that $X =
  X_{0; T,f}$ has fundamental group $\pi$.

  Following the computations in \autocite{torus_fibrations}, we have
  that $\chi(X_0) = 0$ since $\chi(S^1) = 0$ by multiplicativity of
  Euler characteristic.

  Similarly, since both $X_0 \setminus\mathring{\nu}T$ and $T^2
  \times D^2$ have Euler characteristic $0$ (since $T$ is a
  torus-fibered $T^2$-knot), the gluing formula for Euler
  characteristic implies that $\chi(X) = 0$.

  We will show that $\chi(X) > 0$.

  While $\pi$ itself is torsion-free by Lemma~\ref{lem:asph}, its
  abelianization $\pi^\ab = H_1(Y; \Z)$ may in general contain
  torsion. Let $F \cong \Z^{\beta_1(\pi)}$ be the free part of
  $\pi^\ab$ (which has rank $\beta_1(\pi) \geq 1$ by assumption), and
  let $T$ be the torsion subgroup of $\pi^\ab$, so that $\pi^\ab
  \cong F \oplus T$.

  By a routine computation, we have the following homology and
  cohomology groups:

  $$
  H_i(Y) \cong
  \begin{cases}
    \Z, & i = 0, 3\\
    F \oplus T, & i = 1,\\
    F, &i=2\\
    0, & i\geq 4
  \end{cases}, \quad\quad
  H^i(Y;\Z) \cong
  \begin{cases}
    \Z, & i=0,3\\
    F, & i=1,\\
    F \oplus T, &i=2\\
    0, & i\geq 4
  \end{cases}
  $$
  and
  $$
  H_i(X_0) \cong
  \begin{cases}
    \Z, & i=0,4\\
    F \oplus T \oplus \Z, &i=1\\
    F \oplus F \oplus T, &i=2\\
    F \oplus \Z, &i=3\\
    0, &i\geq 5
  \end{cases}, \quad\quad %
  H^i(X_0;\Z) \cong
  \begin{cases}
    \Z, &i=0,4\\
    F \oplus \Z, &i=1\\
    F \oplus F \oplus T, &i=2\\
    F \oplus T \oplus \Z, &i=3\\
    0, & i\geq 5
  \end{cases}.
  $$

  Since $Y$ is a $K(\pi, 1)$ (and is thus a classifying space $B\pi$)
  and $X$ has fundamental group $\pi$, it follows that there is a
  classifying map $\kappa: X \to Y$ inducing an isomorphism on
  $\pi_1$. Since the universal cover $\widetilde{X}$ is simply
  connected, $H^1(\widetilde{X}; \Z) = 0$, so the five-term exact
  sequence in low-degree cohomology associated to the universal cover
  yields an injective map
  $$\kappa^\ast: H^2(Y;\Z) \hookrightarrow H^2(X; \Z).$$

  Restricting $\kappa^\ast$ to the free summand $F \subseteq H^2(Y;
  \Z)$, its image $L = \kappa^\ast(F)$ is a free abelian subgroup of
  $H^2(X;\Z)$ with rank $\beta_1(\pi)$. Moreover, for any $u,v \in
  H^2(Y; \Z)$, their cup product $u \smile v$ lies in $H^4(Y; \Z) =
  0$; by naturality of the cup product, the commutative diagram in
  Figure~\ref{fig:cd2}
  \begin{figure}[!h]
    \centering
    \[
      \begin{tikzcd}[column sep=large]
        H^2(Y; \Z) \times H^2(Y; \Z) \arrow[r, "\smile"] \arrow[d,
        "\kappa^\ast \times \kappa^\ast"', hook]
        \& H^4(Y; \Z) = 0 \arrow[d, "\kappa^\ast"] \\
        H^2(X; \Z) \times H^2(X; \Z) \arrow[r, "\smile"]
        \& H^4(X; \Z) \arrow[r, "{\gen{-, \br{X}}}", "\cong"']
        \& \Z
    \end{tikzcd}\]
    \caption{Cup Product Diagram}\label{fig:cd2}
  \end{figure}
  shows that
  $$Q_X(\kappa^\ast(u), \kappa^\ast(v)) = \gen{\kappa^\ast(u \smile
  v), \br{X}} = 0.$$

  Thus $L$ is isotropic with respect to the intersection form $Q_X$.
  Since $Q_X$ is unimodular and $L$ is isotropic, there is a dual
  free abelian subgroup $L^\ast$ inside of $H^2(X;\Z)$ satisfying
  \mbox{$\rk L^\ast = \rk L = \beta_1(\pi)$} and $L \cap L^\ast = \set{0}$
  (equivalently, the orthogonal complement $L^\perp$ in
    $H^2(X;\Z)/\Tor$ contains $L$ and has rank \mbox{$\rk H^2(X;\Z) - \rk L
  \geq \rk L$}). Hence we obtain that $$\rk H^2(X;\Z) \geq \rk L + \rk
  L^\ast =  2 \beta_1(\pi).$$

  Computing $\chi(X)$ (using Poincar\'{e} duality) we have that:
  $$
  \chi(X) = 2 - 2\beta_1(\pi) + \beta_2(X) \geq 2 - 2\beta_1(\pi) + 2
  \beta_1(\pi) = 2.
  $$

  Hence $\chi(X) \geq 2$. This contradicts the fact that $\chi(X) = 0$.

  Hence there cannot be a torus-fibered $T^2$-knot $T$ in $X_0$ with
  the prescribed surgery.
\end{proof}

Notice that the second half of the proof of Theorem~\ref{thm:main}
shows more generally that any closed, connected, oriented
$4$-manifold $X$ with fundamental group $\pi = \pi_1(Y)$ must satisfy
\mbox{$\beta_2(X) \geq 2\beta_1(\pi)$} and hence $\chi(X) \geq 2$. In
particular, no surgery operation on \mbox{$X_0 = Y \times S^1$} that
preserves the vanishing of the Euler characteristic --- including any
generalized \hbox{Fintushel--Stern} knot surgery along a torus-fibered
$T^2$-knot --- can produce a $4$-manifold with fundamental group $\pi$.

\section{Examples}\label{sec:examples}

We now illustrate Theorem~\ref{thm:main} across several geometric
families of $3$- and $4$-manifolds, beginning with a non-abelian
nilmanifold example and then proceeding to flat (Euclidean
$\mathbb{E}^3$) manifolds and surface bundles over the circle.

\begin{example}\label{ex:heis}
  Let \mbox{$\pi = \pres{x,y,t}{xt=tx,yt=ty, xy=tyx}$} be the
  $3$-dimensional discrete Heisenberg group. This group $\pi$ is the
  fundamental group of the closed, connected, orientable, aspherical
  $3$-dimensional Heisenberg nilmanifold $M$.

  The abelianization is $\pi^\ab \cong \Z^2$, generated by the images
  of $x$ and $y$, so that $\beta_1(\pi) = 2 \geq 1$. By
  Theorem~\ref{thm:main}, no torus surgery along a torus-fibered
  $T^2$-knot in $M \times S^1$ can produce a $4$-manifold with
  fundamental group $\pi$.

  In this example, the manifold $X_0 = M \times S^1$ is the
  Kodaira-Thurston manifold \autocite{Thurston1976}, which is a
  complex, symplectic $4$-manifold that is not K\"{a}hler.
\end{example}

Whereas the Heisenberg manifold $M$ in Example~\ref{ex:heis} is a
non-trivial circle bundle over $T^2$ carrying $\mathrm{Nil}$
geometry, the simplest flat analog is the $3$-torus $T^3$, whose
product with $S^1$ is the $4$-torus $T^4$.

\begin{example}\label{ex:torus}
  Let $Y = T^3$ be the $3$-torus, which is a closed, connected,
  oriented, flat (and hence aspherical) $3$-manifold with \mbox{$\pi =
  \pi_1(T^3) \cong \Z^3$} and $\beta_1(\pi) = 3 \geq 1$. Here \mbox{$X_0 =
  T^3 \times S^1 = T^4$} is the $4$-torus.

  Any coordinate $2$-torus \mbox{$T = T^2 \times \set{\mathrm{pt}} \subset
  T^4$} is a torus-fibered $T^2$-knot, as its exterior \mbox{$T^4 \setminus
  \mathring{\nu}T \cong T^2 \times \Sigma_{1,1}$} is a trivial
  $\Sigma_{1,1}$-bundle over $T^2$. More generally, for any fibered
  knot $K \subseteq T^3$, the torus \mbox{$T = K \times S^1 \subseteq T^3
  \times S^1$} is a torus-fibered $T^2$-knot
  \autocite{torus_fibrations}. By Theorem~\ref{thm:main}, no torus
  surgery on $T^4$ along any of these torus-fibered $T^2$-knots can
  yield a $4$-manifold with fundamental group $\Z^3 = \pi_1(T^3)$.
\end{example}

Both of the preceding examples can be viewed as $T^2$-bundles over
$S^1$. More generally, replacing the torus fiber with a closed
surface $\Sigma_g$ of genus $g \geq 1$ produces a large class of
aspherical $3$-manifolds with positive first Betti number to which
Theorem~\ref{thm:main} applies.

\begin{example}\label{ex:surface_bundles}
  Let $\Sigma_g$ be a closed, connected, oriented surface of genus $g
  \geq 1$, and let \mbox{$Y = \Sigma_g \tilde{\times}_\varphi S^1$} be the
  mapping torus of an orientation-preserving diffeomorphism \mbox{$\varphi:
  \Sigma_g \to \Sigma_g$.} Such a manifold $Y$ can typically be
  obtained as $0$-framed Dehn surgery along a non-trivial fibered
  knot $K$ in $S^3$.

  Because $\Sigma_g$ ($g \geq 1$) and $S^1$ are aspherical, the long
  exact sequence of homotopy groups for the bundle \mbox{$\Sigma_g \to Y
  \onto S^1$} implies that $Y$ is aspherical. Now, depending on the
  \hbox{Nielsen--Thurston} type of $\varphi$, $Y$ carries either $\mathbb{H}^2
  \times \R$, $\mathrm{Solv}$, or $\mathbb{H}^3$ geometry
  \autocite{thurston82}. Furthermore, the bundle projection $Y \onto
  S^1$ induces a surjection \mbox{$\pi = \pi_1(Y) \onto \Z$,} so
  $\beta_1(\pi) \geq 1$. By Theorem~\ref{thm:main}, no torus surgery
  along a torus-fibered $T^2$-knot in $X_0 = Y \times S^1$ can
  produce a $4$-manifold with fundamental group $\pi_1(Y)$.
\end{example}

When the monodromy $\varphi$ in Example~\ref{ex:surface_bundles} is
the identity, $X_0 = Y \times S^1$ splits as a product $\Sigma_g
\times T^2$, where torus-fibered $T^2$-knots arise directly as the $T^2$ fibers.

\begin{example}\label{ex:product_surface}
  As a concrete special case of Example~\ref{ex:surface_bundles} with
  trivial monodromy $\varphi = \id$, let $Y = \Sigma_g \times S^1$
  for $g \geq 1$. Then $Y$ is aspherical with \mbox{$\pi = \pi_1(\Sigma_g)
  \times \Z$} and \mbox{$\beta_1(\pi) = 2g + 1 \geq 3$,} and
  \mbox{$X_0 = Y \times
  S^1 = \Sigma_g \times T^2$.}

  In $X_0 = \Sigma_g \times T^2$, a fiber \mbox{$T = \set{\mathrm{pt}}
  \times T^2$} is a torus-fibered $T^2$-knot whose exterior is the
  product bundle \mbox{$\Sigma_{g,1} \times T^2 \onto T^2$.}
  Theorem~\ref{thm:main} implies that no torus surgery along
  $\set{\mathrm{pt}} \times T^2$ (or any other torus-fibered
  $T^2$-knot in $\Sigma_g \times T^2$) can yield a $4$-manifold with
  fundamental group $\pi_1(\Sigma_g) \times \Z$.
\end{example}

We conclude by showing that the asphericity assumption on $Y$ in
Theorem~\ref{thm:main} is essential: if $Y$ has $\beta_1(Y) > 0$ but
fails to be aspherical, then torus surgery along a torus-fibered
$T^2$-knot in $Y \times S^1$ can indeed realize $\pi_1(Y)$.

\begin{example}[Sharpness of the asphericity hypothesis]\label{ex:sharpness}
  Let $Y = S^1 \times S^2$. Then $Y$ is a closed, connected, oriented
  $3$-manifold with \mbox{$\pi = \pi_1(Y) \cong \Z$} and $\beta_1(\pi) = 1
  \geq 1$, but $Y$ is \emph{not} aspherical since $\pi_2(Y) \cong \Z$.

  In \mbox{$X_0 = Y \times S^1 \cong S^2 \times T^2$,} the $2$-torus \mbox{$T =
  \set{\mathrm{pt}} \times T^2$} has self-intersection zero and exterior
  \[
    X_0 \setminus \mathring{\nu}T \cong D^2 \times T^2 = \Sigma_{0,1}
    \times T^2,
  \]
  so $T$ is a torus-fibered $T^2$-knot in $X_0$. Performing torus
  surgery on $X_0$ along $T$ glues \mbox{$X_0 \setminus \mathring{\nu}T
  \cong D^2 \times T^2$} to $T^2 \times D^2$ along their $T^3$
  boundaries. Choosing an orientation-preserving gluing
  diffeomorphism \mbox{$f: \partial(T^2 \times D^2) \to \partial(D^2 \times
  T^2)$} that exchanges the $\partial D^2$ factor with the $S^1$
  factor of $Y = S^1 \times S^2$ (as in \autocite{torus_fibrations})
  produces \mbox{$X = X_{0; T, f} \cong S^1 \times S^3$,} which has
  fundamental group \mbox{$\pi_1(X) \cong \Z \cong \pi$.} Thus, the
  assumption in Theorem~\ref{thm:main} that $Y$ is aspherical cannot be omitted.
\end{example}

\section{Conclusion}\label{sec:conclusion}

Taken together, Theorem~\ref{thm:main} and the examples of
Section~\ref{sec:examples} highlight a fundamental rigidity in the
topology of $4$-manifolds obtained by surgery on $X_0 = Y \times
S^1$. While the fibration extension results of
\autocites{meyerthesis,torus_fibrations} demonstrate that torus
surgery along a torus-fibered $T^2$-knot $T \subseteq X_0$ always
yields a closed $4$-manifold $X = X_{0; T, f}$ that fibers over
$S^1$, the resulting manifold cannot realize the $3$-manifold group
$\pi = \pi_1(Y)$ whenever $Y$ is closed, connected, oriented, and
aspherical with $\beta_1(Y) > 0$. As Example~\ref{ex:sharpness}
shows, this obstruction depends crucially on the asphericity of $Y$,
which allows non-trivial second cohomology classes in $H^2(Y; \Z)$ to
pull back via the classifying map $\kappa: X \to Y$ to an isotropic
sublattice of $H^2(X; \Z)$.

We close by noting two natural extensions of this discussion:
\begin{enumerate}
  \item \emph{The case $\beta_1(Y) = 0$:} If $Y$ is a closed,
    connected, oriented $3$-manifold with $\beta_1(Y) = 0$ (i.e., a
    rational homology $3$-sphere), then $H_1(Y; \Z) = \pi^\ab$ is
    finite. On the other hand, by \autocite{torus_fibrations}, any
    torus surgery $X = X_{0; T, f}$ along a torus-fibered $T^2$-knot
    $T \subseteq X_0$ fibers over $S^1$ and therefore admits a
    surjection $\pi_1(X) \onto \Z$, forcing $\beta_1(X) \geq 1$.
    Consequently, $\pi_1(X_{0; T, f})$ cannot be isomorphic to
    $\pi_1(Y)$ when $\beta_1(Y) = 0$ either. Combining this
    observation with Theorem~\ref{thm:main} shows that for
    \emph{every} closed, connected, oriented, aspherical $3$-manifold
    $Y$ (regardless of $\beta_1(Y)$), no torus surgery on $Y \times
    S^1$ along a torus-fibered $T^2$-knot can yield a $4$-manifold
    with fundamental group $\pi_1(Y)$.
  \item \emph{Arbitrary Euler-characteristic-preserving surgeries:}
    Because the proof of Theorem~\ref{thm:main} uses only the fact
    that $\chi(X) = 0$ together with the classifying map $X \to Y$,
    the same cohomological obstruction applies verbatim to torus
    surgery along \emph{any} self-intersection zero $2$-torus $T
    \subseteq Y \times S^1$ (not necessarily torus-fibered), as well as
    to any generalized \hbox{Fintushel--Stern} knot surgery along a
    self-intersection zero $2$-torus in $Y \times S^1$, whenever $Y$
    is closed, connected, oriented, and aspherical with $\beta_1(Y) > 0$.
\end{enumerate}

\section*{Acknowledgments}
The author would like to thank Susan Hermiller for suggesting
examination of the Heisenberg group in Example~\ref{ex:heis}. The
author would also like to thank Rom\'{a}n Aranda for helpful
correspondence while completing this manuscript.

\subsection*{AI Disclosure:}
During the preparation of this work, the author utilized Gemini 4
Argon to assist with preliminary literature identification and
exploratory computations, alongside Grammarly for grammar correction
and stylistic refinement. After using these tools, the author
independently reviewed, cross-referenced, and verified all cited
literature; independently executed and verified all mathematical
calculations; and thoroughly edited the manuscript. The results,
proofs, and conclusions are entirely the author's own, and the author
assumes full responsibility for the integrity and content of this publication.

\printbibliography

@article{FS98,
  title = {Knots, links, and 4-manifolds},
  volume = {134},
  issn = {1432-1297},
  url = {http://dx.doi.org/10.1007/s002220050268},
  doi = {10.1007/s002220050268},
  number = {2},
  journal = {Inventiones Mathematicae},
  publisher = {Springer Science and Business Media LLC},
  author = {Fintushel, Ronald and Stern, Ronald J.},
  year = {1998},
  month = oct,
  pages = {363--400}
}

@article{Donaldson87,
  title = {Irrationality and the $h$-cobordism conjecture},
  volume = {26},
  issn = {0022-040X},
  url = {http://dx.doi.org/10.4310/jdg/1214441179},
  doi = {10.4310/jdg/1214441179},
  number = {1},
  journal = {Journal of Differential Geometry},
  publisher = {International Press of Boston},
  author = {Donaldson, Simon K.},
  year = {1987},
  month = jan
}

@inbook{Dolgachev2010,
  title = {Algebraic Surfaces with $q = p_g = 0$},
  isbn = {9783642110870},
  url = {http://dx.doi.org/10.1007/978-3-642-11087-0_3},
  doi = {10.1007/978-3-642-11087-0_3},
  booktitle = {Algebraic Surfaces},
  publisher = {Springer Berlin Heidelberg},
  author = {Dolgachev, Igor},
  year = {2010},
  pages = {97--215}
}

@article{Zentner,
  title = {A note on logarithmic transformations on the {Hopf} surface},
  volume = {342},
  issn = {1778-3569},
  url = {http://dx.doi.org/10.1016/j.crma.2005.10.030},
  doi = {10.1016/j.crma.2005.10.030},
  number = {1},
  journal = {Comptes Rendus. Math\'{e}matique},
  publisher = {Cellule MathDoc/Centre Mersenne},
  author = {Zentner, Raphael},
  year = {2005},
  month = nov,
  pages = {37--42}
}

@article{Ni17,
  title = {{Fintushel--Stern} knot surgery in torus bundles},
  volume = {10},
  issn = {1753-8424},
  url = {http://dx.doi.org/10.1112/topo.12002},
  doi = {10.1112/topo.12002},
  number = {1},
  journal = {Journal of Topology},
  publisher = {Wiley},
  author = {Ni, Yi},
  year = {2017},
  month = jan,
  pages = {164--177}
}

@article{Kim19,
  title = {Knotted surfaces in $4$-manifolds by knot surgery and stabilization},
  volume = {371},
  issn = {1088-6850},
  url = {http://dx.doi.org/10.1090/tran/7450},
  doi = {10.1090/tran/7450},
  number = {12},
  journal = {Transactions of the American Mathematical Society},
  publisher = {American Mathematical Society (AMS)},
  author = {Kim, Hee Jung},
  year = {2019},
  month = feb,
  pages = {8405--8427}
}

@article{Akbulut2002,
  title = {Variations on {Fintushel-Stern} Knot Surgery on $4$-manifolds},
  volume = {26},
  number = {1},
  journal = {Turkish Journal of Mathematics},
  author = {Akbulut, Selman},
  year = {2002},
  url = {https://journals.tubitak.gov.tr/math/vol26/iss1/6},
  pages = {81--92}
}

@article{Yun2008,
  title = {Twisted fiber sums of {Fintushel-Stern's} knot surgery 4-manifolds},
  volume = {360},
  issn = {1088-6850},
  url = {http://dx.doi.org/10.1090/S0002-9947-08-04623-0},
  doi = {10.1090/s0002-9947-08-04623-0},
  number = {11},
  journal = {Transactions of the American Mathematical Society},
  publisher = {American Mathematical Society (AMS)},
  author = {Yun, Ki-Heon},
  year = {2008},
  month = may,
  pages = {5853--5868}
}

@article{AEMS2008,
  title = {A note on {Stein} fillings of contact manifolds},
  volume = {15},
  issn = {1945-001X},
  url = {http://dx.doi.org/10.4310/MRL.2008.v15.n6.a5},
  doi = {10.4310/mrl.2008.v15.n6.a5},
  number = {6},
  journal = {Mathematical Research Letters},
  publisher = {International Press of Boston},
  author = {Akhmedov, Anar and Etnyre, John B. and Mark, Thomas E.
  and Smith, Ivan},
  year = {2008},
  pages = {1127--1132}
}

@misc{torus_fibrations,
  doi = {10.48550/ARXIV.2512.17595},
  url = {https://arxiv.org/abs/2512.17595},
  eprint = {2512.17595},
  author = {Meyer, Nicholas},
  title = {Extending Fibrations of the $3$-Torus and Applications to
  Torus Surgery in $4$-Manifolds},
  publisher = {arXiv},
  year = {2025},
  copyright = {arXiv.org perpetual, non-exclusive license},
  note = {\href{https://arxiv.org/abs/2512.17595}{arXiv:2512.17595}} 
}

@book{hatcher,
  address = {New York},
  title = {Algebraic topology},
  isbn = {978-0-521-79540-1},
  publisher = {Cambridge University Press},
  author = {Hatcher, Allen},
  year = {2001}
}

@article{thurston82,
  author = {Thurston, William P.},
  title = {{Three dimensional manifolds, Kleinian groups and
  hyperbolic geometry}},
  volume = {6},
  journal = {Bulletin (New Series) of the American Mathematical Society},
  number = {3},
  publisher = {American Mathematical Society},
  pages = {357--381},
  year = {1982}
}

@phdthesis{meyerthesis,
  title = {Torus Surgery, Fibrations, Multisections, and Spun 4-Manifolds},
  author = {Meyer, Nicholas Paul},
  year = {2024},
  month = jul,
  school = {University of Nebraska --- Lincoln},
  url = {https://digitalcommons.unl.edu/dissunl/1150/}
}

@article{Scott73,
  title = {Compact Submanifolds of 3-Manifolds},
  volume = {s2-7},
  issn = {0024-6107},
  url = {http://dx.doi.org/10.1112/jlms/s2-7.2.246},
  doi = {10.1112/jlms/s2-7.2.246},
  number = {2},
  journal = {Journal of the London Mathematical Society},
  publisher = {Wiley},
  author = {Scott, G. P.},
  year = {1973},
  month = nov,
  pages = {246--250}
}

@article{Papa57,
  title = {On {Dehn's} Lemma and the Asphericity of Knots},
  volume = {66},
  issn = {0003-486X},
  url = {http://dx.doi.org/10.2307/1970113},
  doi = {10.2307/1970113},
  number = {1},
  journal = {The Annals of Mathematics},
  publisher = {JSTOR},
  author = {Papakyriakopoulos, C. D.},
  year = {1957},
  month = jul,
  pages = {1--32}
}

@article{Whitehead1,
  title = {Combinatorial homotopy. {I}},
  volume = {55},
  issn = {0273-0979},
  url = {http://dx.doi.org/10.1090/S0002-9904-1949-09175-9},
  doi = {10.1090/s0002-9904-1949-09175-9},
  number = {3},
  journal = {Bulletin of the American Mathematical Society},
  publisher = {American Mathematical Society (AMS)},
  author = {Whitehead, J. H. C.},
  year = {1949},
  pages = {213--245}
}

@article{Whitehead2,
  title = {Combinatorial homotopy. {II}},
  volume = {55},
  issn = {0273-0979},
  url = {http://dx.doi.org/10.1090/S0002-9904-1949-09213-3},
  doi = {10.1090/s0002-9904-1949-09213-3},
  number = {5},
  journal = {Bulletin of the American Mathematical Society},
  publisher = {American Mathematical Society (AMS)},
  author = {Whitehead, J. H. C.},
  year = {1949},
  pages = {453--496}
}

@article{Thurston1976,
  title = {Shorter Notes: Some Simple Examples of Symplectic Manifolds},
  volume = {55},
  issn = {0002-9939},
  url = {http://dx.doi.org/10.2307/2041749},
  doi = {10.2307/2041749},
  number = {2},
  journal = {Proceedings of the American Mathematical Society},
  publisher = {JSTOR},
  author = {Thurston, William P.},
  year = {1976},
  month = mar,
  pages = {467--468}
}

\end{document}